\documentclass[final,1p,times]{elsarticle}

\usepackage{amssymb}
 \usepackage{amsthm}
\usepackage{amscd}
\usepackage{bm}
\usepackage{amsmath}
\usepackage{amsfonts}
\usepackage{amssymb}
\usepackage{graphicx}
\newtheorem{theorem}{Theorem}

\newtheorem{lemma}[theorem]{Lemma}

\usepackage{mathrsfs}
\usepackage{titletoc}

\newcommand{\ra}{\rightarrow}
\newcommand{\p}{\partial}
\newcommand{\f}{\frac}

\newcommand{\be}{\begin{equation}}
\renewcommand{\ra}{\rightarrow}
\newcommand{\ee}{\end{equation}}
\newcommand{\bea}{\begin{eqnarray}}
\newcommand{\eea}{\end{eqnarray}}
\newcommand{\bna}{\begin{eqnarray*}}
\newcommand{\ena}{\end{eqnarray*}}

\renewcommand{\le}{\left}
\newcommand{\ri}{\right}

\journal{***}

\begin{document}

\begin{frontmatter}

\title{Prescribing Gaussian curvature on closed Riemann surface with conical singularity in the negative case}

 \author{Yunyan Yang}
 \ead{yunyanyang@ruc.edu.cn}
 \author{Xiaobao Zhu}
 \ead{zhuxiaobao@ruc.edu.cn}

\address{ Department of Mathematics,
Renmin University of China, Beijing 100872, P. R. China}

\begin{abstract}
  The problem of prescribing Gaussian curvature on Riemann surface with conical singularity is considered.
  Let $(\Sigma,\bm\beta)$ be a closed Riemann surface with a divisor $\bm\beta$, and $K_\lambda=K+\lambda$,
  where $K:\Sigma\ra\mathbb{R}$ is a H\"older continuous function satisfying $\max_\Sigma  K= 0$, $K\not\equiv 0$,
  and $\lambda\in\mathbb{R}$.
  If the Euler characteristic $\chi(\Sigma,\bm\beta)$ is negative, then by a variational
  method, it is proved that
  there exists a constant $\lambda^\ast>0$ such that
   for any $\lambda\leq 0$, there is a unique conformal metric with the Gaussian curvature $K_\lambda$;  for any $\lambda$, $0<\lambda<\lambda^\ast$, there are at least two conformal metrics having $K_\lambda$ its Gaussian curvature;
   for $\lambda=\lambda^\ast$,
  there is at least one conformal metric with the Gaussian curvature $K_{\lambda^\ast}$;
   for any $\lambda>\lambda^\ast$, there is no certain conformal metric having $K_{\lambda}$
  its Gaussian curvature. This result is an analog of that of Ding and Liu \cite{Ding-Liu}, partly resembles
  that of Borer, Galimberti and Struwe \cite{B-G-Stru}, and generalizes that of Troyanov \cite{Troyanov} in the negative case.
\end{abstract}


\end{frontmatter}

\titlecontents{section}[0mm]
                       {\vspace{.2\baselineskip}}%\bfseries}
                       {\thecontentslabel~\hspace{.5em}}
                        {}
                        {\dotfill\contentspage[{\makebox[0pt][r]{\thecontentspage}}]}
\titlecontents{subsection}[3mm]
                       {\vspace{.2\baselineskip}}%\bfseries}
                       {\thecontentslabel~\hspace{.5em}}
                        {}
                       {\dotfill\contentspage[{\makebox[0pt][r]{\thecontentspage}}]}

\setcounter{tocdepth}{2}
%\tableofcontents

%\setcounter{tocdepth}{1}
%\tableofcontents

%\tableofcontents
\section{Introduction}
The problem of prescribing Gaussian curvature on smooth Riemann surfaces has been well understood
\cite{KW1}. Let $(\Sigma,g)$ be a closed smooth Riemann surface, $\chi(\Sigma)$ be its topological Euler characteristic,
and $\kappa:\Sigma\ra\mathbb{R}$ be its Gaussian curvature. If $\bar{g}=e^{2u}g$ with a smooth function $u$, then
the Gaussian curvature of $(\Sigma,\bar{g})$ satisfies
 $\bar{\kappa}=e^{-2u}({\kappa}+\Delta_gu)$, where $\Delta_g$ denotes the Laplacce-Beltrami
operator with respect to the metric $g$. A natural question is whether for any smooth function $K:\Sigma\ra\mathbb{R}$,
there is a smooth function $u$ such that the metric $e^{2u}g$ has $K$ its Gaussian curvature.
Clearly this is equivalent to solving the elliptic equation
\be\label{eqa}\Delta_gu+\kappa-Ke^{2u}=0\quad{\rm on}\quad \Sigma.\ee
The Gauss-Bonnet formula leads to
$$\int_\Sigma Ke^{2u}dv_g=\int_\Sigma \kappa dv_g=2\pi\chi(\Sigma).$$
Note that the solvability of (\ref{eqa}) is closely related to the sign of $\chi(\Sigma)$.
If $\chi(\Sigma)>0$, then $\Sigma$ is either the projective space $\mathbb{RP}^2$
or the $2$-sphere $\mathbb{S}^2$. In the case of $\mathbb{RP}^2$,
it was shown by Moser \cite{Moser} that the equation (\ref{eqa}) has a solution $u$, provided that $K\in C^\infty(\mathbb{S}^2)$
satisfies $\sup_\Sigma K>0$ and $K(p)=K(-p)$ for all $p\in \mathbb{S}^2$. While the problem on $\mathbb{S}^2$ is much more
complicated and known as the Nirenberg problem, see for examples \cite{KW1,C-Y1,C-Y2,C-Liu,Chen-Ding}. If $\chi(\Sigma)=0$, the problem
has been completely solved by Kazdan-Warner \cite{KW1}. While if $\chi(\Sigma)<0$, the problem was studied by Kazdan and Warner \cite{KW1}
via the method of upper and lower solutions. They proved that if $K\leq 0$ and $K\not\equiv 0$, then (\ref{eqa}) has a unique solution. Later, Ding and Liu
\cite{Ding-Liu} considered the case that $K$ changes sign. Precisely, replacing $K$ by $K+\lambda$ in (\ref{eqa}) with $K\leq 0$,
$K\not\equiv 0$, and $\lambda\in\mathbb{R}$, they obtained the following conclusion by using
a method of upper and lower solutions and a variational method: there exists a $\lambda^\ast>0$ such that  if $\lambda\leq 0$, then
(\ref{eqa}) has a unique solution;  if $0<\lambda<\lambda^\ast$, then (\ref{eqa}) has at least two solutions;
if $\lambda=\lambda^\ast$, then (\ref{eqa}) has at least one solution;  if $\lambda>\lambda^\ast$, then (\ref{eqa}) has no solution.
Recently, using a monotonicity technique due to Struwe \cite{Struwe1,Struwe2}, Borer, Galimberti,
and Struwe \cite{B-G-Stru} partly reproved the above results and obtained additional
estimates for certain sequence of solutions that allow to characterize their bubbling behavior. Further analysis
in this direction has been done
by Galimberti \cite{Galim}, del Pino and Rom\'an \cite{del}.

The problem of prescribing Gaussian curvature can also be proposed on surfaces with conical singularities.
Let $\Sigma$ be a closed Riemann surface, $p_1,\cdots,p_\ell$ be points
of $\Sigma$ and $\theta_1,\cdots,\theta_\ell$ be positive numbers. Denote
$$\bar\chi=2\pi \chi(\Sigma)+\sum_{i=1}^\ell(\theta_i-2\pi).$$ Then it was proved by Troyanov \cite{Troyanov}
that
if $0<\bar\chi<\min\{4\pi,2\theta_1,\cdots,2\theta_\ell\}$,
then any smooth function on $\Sigma$, which is positive at some point is the Gaussian curvature of a conformal metric having
at $p_i$ a conical singularity of angle $\theta_i$; if $\bar\chi=0$, then a smooth nonconstant
function $K:\Sigma\ra\mathbb{R}$ is the Gaussian curvature of a conformal metric having
at $p_i$ a conical singularity of angle $\theta_i$ if and only if $\int_\Sigma Kd\mu<0$, where $d\mu$ is the area element of
the original singular metric; if $\bar\chi<0$, then any smooth negative function on $\Sigma$ is the Gaussian curvature
of a unique conformal metric having at $p_i$ a conical singularity of angle $\theta_i$. As in the smooth Riemann surface case,
the prescribing Gaussian curvature problem on the $2$-sphere with conical singularity is most delicate.
The case $\ell=2$ was studied by Chen and Li \cite{CL1,CL2}.
While the case $\ell\geq 3$ was considered by Eremenko
\cite{Eremenko}, Malchiodi and Ruiz \cite{Malchiodi}, Chen and Lin \cite{ChLin},  Marchis
and L\'opez-Soriano \cite{LS}, and others.

In this paper, we focus on the negative case, namely $\bar\chi<0$. Precisely we shall prove an analog of the result of
Ding and Liu \cite{Ding-Liu}, and thereby part of results of Borer, Galimberti,
and Struwe \cite{B-G-Stru}. Though we still use the variational method, which had been employed by Ding and Liu, we have to overcome
difficulties in the presence of conical singularities. In particular, we have to establish the strong maximum principle, which
is essential for the method of upper and lower solutions in our setting.

The remaining part of this paper is organized as follows: In Section 2, we give some notations for surfaces with conical singularities
and state our main results;
In Section 3, the maximum principle for the Laplace-Beltrami operator and the Palais-Smale condition
for certain functional are discussed; In Section 4, following the lines of
\cite{Ding-Liu,B-G-Stru}, we prove our main theorem.

\section{Notations and main results}
Let us briefly recall some geometric concepts from Troyanov \cite{Troyanov}.
In general, a closed Riemann surface $\Sigma$ is defined to be a topological space with an atlas
$\{\phi_i:U_i\ra\mathbb{C}\}$, where if $U_i\cap U_j\not=\varnothing$, then the coordinate transformation
$\phi_i\circ \phi_j^{-1}$ is conformal, i.e., holomorphic or anti-holomorphic. Two such atlases define the same
structure on $\Sigma$ if their union is still such an atlas. A conformal Riemannian metric is defined by
$g=\rho(z)|dz|^2$ locally, where $z$ is a coordinate on $\Sigma$ and $\rho$ is a positive measurable function.
A divisor on a Riemann surface is a formal sum
${\bm\beta}=\sum_{i=1}^\ell\beta_ip_i$,
where $p_i\in \Sigma$ and $\beta_i>-1$, $i=1,\cdots,\ell$. The set ${\rm supp}\,\bm\beta=\{p_1,\cdots,p_\ell\}$ is the support of $\bm\beta$,
and the number $|\bm\beta|=\sum_{i=1}^\ell\beta_i$ is the degree of the divisor.
A conformal metric $g$ on $\Sigma$ is said to represent the divisor $\bm\beta$ if $g\in C^2(\Sigma\setminus{\rm supp}\,\bm\beta)$ verifying that
if $z_i$ is a coordinate defined in a neighborhood $U_i$ of $p_i$, then there is some $u_i\in C^2(U_i\setminus\{p_i\})\cap
C^0(U_i\setminus\{p_i\})$ such that
\be\label{metric}g=e^{2u_i}|z_i-z_i(p_i)|^{2\beta_i}|dz_i|^{2}.\ee
Under the circumstances, $g$ is said to have a conical singularity of order $\beta_i$ or angle $\theta_i=2\pi(\beta_i+1)$ at $p_i$,
$i=1,\cdots,\ell$.
The Euler characteristic of $(\Sigma,\bm\beta)$ is defined by
$$\chi(\Sigma,\bm\beta)=\chi(\Sigma)+|\bm\beta|,$$
where $\chi(\Sigma)$ is the topological Euler characteristic of $\Sigma$, and $|\bm\beta|=\sum_{i=1}^\ell \beta_i$ is the degree of $\bm\beta$.
Let $\kappa: \Sigma\setminus {\rm supp}\,
\bm\beta\ra \mathbb{R}$ be the Gaussian curvature of $g$. If $\kappa$ can be extended to a H\"older continuous function on $\Sigma$, then
it was shown by Troyanov \cite{Troyanov} that
a Gauss-Bonnet formula holds:
\be\label{Gauss-Bonnet}\int_\Sigma \kappa dv_g=2\pi\chi(\Sigma,\bm\beta),\ee
where $dv_g$ denotes the Riemannian volume element with respect to the conical metric $g$.

Let $(\Sigma,\bm\beta)$ be a closed Riemann surface with a divisor $\bm\beta=\sum_{i=1}^\ell \beta_ip_i$, and the metric
$g$ represents $\bm\beta$ with $\beta_i>-1$, $i=1,\cdots,\ell$.
It follows from (\ref{metric}) that there exists a smooth Riemannian metric $g_0$ such that
\be\label{g0}g=\rho g_0,\ee where
$\rho>0$ on $\Sigma$, $\rho\in C^2_{\rm loc}(\Sigma\setminus{\rm supp}\,\bm\beta)$,
and $\rho\in L^{r}(\Sigma)$ for some $r>1$.
 Let $W^{1,2}(\Sigma,g)$ be the completion of $C^\infty(\Sigma)$ under the norm
$$\|u\|_{W^{1,2}(\Sigma,g)}=\le(\int_\Sigma(|\nabla_gu|^2+u^2)dv_g\ri)^{1/2},$$
where $\nabla_g$ denotes the gradient operator with respect to the metric $g$.
It was observed by Troyanov \cite{Troyanov} that $W^{1,2}(\Sigma,g)=W^{1,2}(\Sigma,g_0)$.
As a consequence, by the Sobolev embedding theorem for
smooth Riemann surface $(\Sigma,g_0)$ and the H\"older inequality, one has
\be\label{Sobolev}W^{1,2}(\Sigma,g)\hookrightarrow L^p(\Sigma,g),\quad \forall p>1.\ee

We now state the following:
\begin{theorem}\label{mainTheorem}
Let $(\Sigma,\bm\beta)$ be a closed Riemann surface with a divisor $\bm\beta=\sum_{i=1}^\ell\beta_ip_i$.
 Suppose that the Euler characteristic
$\chi(\Sigma,\bm\beta)<0$, $K: \Sigma\ra\mathbb{R}$ is a H\"older continuous function,
$\max_{\Sigma}K= 0$ and $K\not\equiv 0$.
Let $K_\lambda=K+\lambda$, $\lambda\in\mathbb{R}$. Assume that a conformal metric $g$ represents $\bm\beta$. Let
$\kappa:\Sigma\setminus{\rm supp}\,\bm\beta \ra\mathbb{R}$ be the Gaussian curvature of $g$, and $\kappa$
can be extended to a H\"older continuous function on $\Sigma$.
Then there exists a constant $\lambda^\ast>0$ such that
$(i)$ when $\lambda\leq 0$, there exists a unique conformal metric on $\Sigma$ with Gaussian curvature $K_\lambda$, representing
the divisor $\bm\beta$;
$(ii)$ when $0<\lambda<\lambda^\ast$, there exist at least two conformal metrics on $\Sigma$ with the same Gaussian curvature $K_\lambda$,
representing the divisor $\beta$;
$(iii)$ when $\lambda=\lambda^\ast$, there exists at least one conformal metric on $\Sigma$ with Gaussian curvature $K_{\lambda^\ast}$, representing
the divisor $\bm\beta$;
$(iv)$ when $\lambda>\lambda^\ast$, there is no function $u\in W^{1,2}(\Sigma,g)\cap C^2(\Sigma\setminus{\rm supp}\,\bm\beta)\cap
C^0(\Sigma)$ such that
 $e^{2u}g$ has the Gaussian curvature $K_{\lambda}$.
\end{theorem}

Since the metric $g$ has the Gaussian curvature $\kappa$, and the metric $g_\lambda=e^{2u}g$ has the Gaussian curvature $K_\lambda=K+\lambda$.
A standard calculation shows
\be\label{eqn-1}\Delta_gu+\kappa-K_\lambda e^{2u}=0\quad{\rm on}\quad \Sigma\setminus{\rm supp}\,\bm\beta.\ee
 Note that if $u\in W^{1,2}(\Sigma,g)$ is a distributional solution of
the equation
\be\label{eqn}\Delta_gu+\kappa-K_\lambda e^{2u}=0 \quad{\rm on}\quad \Sigma,\ee
we have by elliptic estimates $u\in C^2(\Sigma\setminus{\rm supp}\,\bm\beta)\cap C^0(\Sigma)$, and thus (\ref{eqn-1})
holds. Hence, in order to prove Theorem \ref{mainTheorem}, it suffices to show  the following:
\begin{theorem}\label{thm2}
Under the same assumptions as in Theorem \ref{mainTheorem}, there exists a $\lambda^\ast>0$ such that $(i)$
if $\lambda\leq 0$, then (\ref{eqn}) has a unique distributional solution; $(ii)$ if $0<\lambda<\lambda^\ast$, then (\ref{eqn}) has at least
two distributional solutions; $(iii)$ if $\lambda=\lambda^\ast$, then (\ref{eqn}) has at least one distributional solution;
$(iv)$ if $\lambda>\lambda^\ast$, then (\ref{eqn}) has no distributional solution.
\end{theorem}

For the proof of Theorem \ref{thm2}, we  follow closely  Ding and Liu \cite{Ding-Liu} by employing a variational method.
In particular we use the upper and lower solutions principle and the strong maximum principle.
In the remaining part of this paper, $(\Sigma,g)$ will always denote a conical singular Riemann surface given in Theorem \ref{mainTheorem};
we do  not distinguish sequence and subsequence;
moreover we often denote various constants by the same $C$, even in the same line.

%For any $\lambda\in\mathbb{R}$, we define a functional $E_\lambda:W^{1,2}(\Sigma,g)\ra \mathbb{R}$ by
%\be\label{functional}E_\lambda(u)=\int_\Sigma(|\nabla_gu|^2+2ku-K_\lambda e^{2u})dv_g.\ee

\section{Preliminary analysis}

In this section, we prove maximum principle, Palais-Smale condition, upper and lower solutions principle,
which will be used later. Compared with the smooth Riemann surface case, all the above mentioned things need to be
re-established since the metric $g$ has conical singularity.

\subsection{Maximum principle}

We first have a weak maximum principle by integration by parts, namely

\begin{lemma}[Weak maximum principle]\label{maximum principle}
 For any constant $c>0$, if $u\in W^{1,2}(\Sigma,g)\cap C^0(\Sigma)$ satisfies $\Delta_gu+cu\geq 0$ in the distributional sense, then $u\geq 0$ on $\Sigma$.
\end{lemma}
\noindent{\it Proof}. Denote $u^-=\min\{u,0\}$. Testing the equation $\Delta_gu+cu\geq 0$ by $u^-$, one has
$$\int_\Sigma(|\nabla_gu^-|^2+c{u^-}^2)dv_g\leq 0.$$
This leads to $u^-\equiv 0$ on $\Sigma$. $\hfill\Box$\\

Moreover, using the Moser iteration (see for example Theorems 8.17
and 8.18 in \cite{GT}),  we obtain the following strong maximum principle.

\begin{lemma}[Strong maximum principle]\label{strong-maximum}
 Let $u\in W^{1,2}(\Sigma,g)\cap C^0(\Sigma)$ satisfy that $u\geq 0$ on $\Sigma$, and that
 for some positive constant $c$, $\Delta_gu+cu\geq 0$ in the distributional sense. If there exists a point $x_0\in\Sigma$ such that
 $u(x_0)=0$, then there holds $u\equiv 0$ on $\Sigma$.
\end{lemma}

\noindent{\it Proof}.
{\bf Step 1}. {\it If $v\in W^{1,2}(\Sigma,g)\cap C^0(\Sigma)$ satisfies $v\geq 0$ on $\Sigma$, and
\be\label{low-}\Delta_gv-cv\leq 0\ee
in the distributional sense, where $c$ is a positive constant,
then there exists some constant $C$ depending only on $(\Sigma,g)$ such that
\be\label{Moser-iter}\|v\|_{L^\infty(\Sigma)}\leq C\|v\|_{L^2(\Sigma,g)}.\ee}

Now we use the Moser iteration to prove (\ref{Moser-iter}). For any $p\geq 2$, testing (\ref{low-}) by $v^{p-1}$ and integrating by parts, we have
$$\int_\Sigma |\nabla_gv^{\f{p}{2}}|^2dv_g\leq \f{cp^2}{4(p-1)}\int_\Sigma v^pdv_g.$$
Hence
$\|v^{\f{p}{2}}\|_{W^{1,2}(\Sigma,g)}\leq Cp\|v^{\f{p}{2}}\|_{L^{2}(\Sigma,g)}$ for some constant $C$.
Then the Sobolev embedding (\ref{Sobolev}) leads to
$\|v^{\f{p}{2}}\|_{L^4(\Sigma,g)}\leq Cp\|v^{\f{p}{2}}\|_{L^{2}(\Sigma,g)}$, which is equivalent to
$\|v\|_{L^{2p}(\Sigma,g)}\leq C^{\f{2}{p}}p^{\f{2}{p}}\|v\|_{L^p(\Sigma,g)}$. Taking $p=p_k=2^k$, $k=1,2,\cdots$, we have
\be\label{it}\|v\|_{L^{p_{k+1}}(\Sigma,g)}\leq C^{\f{2}{p_k}}p_k^{\f{2}{p_k}}\|v\|_{L^{p_k}(\Sigma,g)}\leq
C^{\sum_{j=1}^k2^{1-j}}2^{\sum_{j=1}^k2^{1-j}j}\|v\|_{L^2(\Sigma,g)}\leq C\|v\|_{L^2(\Sigma,g)}.\ee
Letting $k\ra\infty$ in (\ref{it}), we conclude (\ref{Moser-iter}).\\

\noindent{\bf Step 2}. {\it Let $u\in W^{1,2}(\Sigma,g)\cap C^0(\Sigma)$ be a nonnegative distributional solution of
\be\label{dis}\Delta_g u+cu\geq 0,\ee
where $c$ is a positive constant.
Then there exists some constant $C$ such that
\be\label{Harnack}\|u\|_{L^2(\Sigma,g)}\leq C\inf_{\Sigma}u.\ee}

Without loss of generality, we assume  $u\geq \epsilon>0$, otherwise we can replace $u$ by $u+\epsilon$.
We {\it claim} that that $u^{-1}$ is a distributional solution of $\Delta_gu^{-1}-cu^{-1}\leq 0$. To see it, we recall that
$g=\rho g_0$, where $\rho:\Sigma\ra\mathbb{R}$ is a positive function, $\rho\in L^q(\Sigma)$ for some $q>1$, and $g_0$ is a smooth Riemannian metric. Then for any
$\phi\in W^{1,2}(\Sigma,g_0)$ with $\phi\geq 0$, we calculate
\bna\int_\Sigma\le(\nabla_{g}u^{-1}\nabla_{g}\phi-c u^{-1}\phi\ri)dv_{g}&=&
\int_\Sigma\le(\nabla_{g_0}u^{-1}\nabla_{g_0}\phi-c\rho u^{-1}\phi\ri)dv_{g_0}\\
&=&-\int_\Sigma\le(\nabla_{g_0}u\nabla_{g_0}(\phi u^{-2})+2\phi u^{-3}|\nabla_{g_0}u|^2+
c\rho u(\phi u^{-2})\ri)dv_{g_0}\\
&\leq& -\int_\Sigma\le(\nabla_{g_0}u\nabla_{g_0}(\phi u^{-2})+c\rho u(\phi u^{-2})\ri)dv_{g_0}\\
&=&-\int_\Sigma\le(\nabla_{g}u\nabla_{g}(\phi u^{-2})+c u(\phi u^{-2})\ri)dv_{g}.
\ena
This together with (\ref{dis}) confirms our claim.
Now we have by Step 1,
$$\sup_\Sigma u^{-1}\leq C\|u^{-1}\|_{L^2(\Sigma,g)},$$
which leads to
\bna
\inf_\Sigma u&\geq&C\le(\int_\Sigma u^{-2}dv_g\ri)^{-\f{1}{2}}\\
&=&C\le(\int_\Sigma u^{-2}dv_g\int_\Sigma u^2dv_g\ri)^{-1/2}\le(\int_\Sigma u^2dv_g\ri)^{1/2}.
\ena
Thus, to prove (\ref{Harnack}), it suffices to show there exists some constant $C$ such that
\be\label{bound}\int_\Sigma u^{-2}dv_g\int_\Sigma u^2dv_g\leq C.\ee
Let $w=\log u-\gamma$, where $\gamma=\f{1}{{\rm Vol}_g(\Sigma)}\int_\Sigma \log u\,dv_g$. We shall prove that
\be\label{exp}\int_\Sigma e^{2|w|}dv_g\leq C,\ee
which implies
$$\int_\Sigma e^{2(\gamma-\log u)}dv_g\leq C,\quad \int_\Sigma e^{2(\log u-\gamma)}dv_g\leq C.$$
This immediately leads to (\ref{bound}).

We are only left to prove (\ref{exp}). Testing the equation (\ref{dis}) by $u^{-1}$, we have
$$\int_\Sigma(\nabla_gu^{-1}\nabla_gu+c)dv_g\geq 0.$$
It follows that
\be\label{nb-w}\int_\Sigma|\nabla_gw|^2dv_g\leq C.\ee
Note that $\int_\Sigma wdv_g=0$. In view of (\ref{nb-w}), we conclude from the Poincar\'e inequality that
 \be\label{w12}\|w\|_{W^{1,2}(\Sigma,g)}\leq C.\ee
 Recall that the metric $g$ represents the divisor $\bm\beta=\sum_{i=1}^\ell \beta_ip_i$ with $\beta_i>-1$, $i=1,\cdots,\ell$.
 Denote $b=\min\{1,1+\beta_1,\cdots,1+\beta_\ell\}$.
 Then the Trudinger-Moser inequality for surfaces with
 conical singularities \cite{Troyanov} together with (\ref{w12}) implies that
\bea\nonumber
\int_\Sigma e^{2|w|}dv_g&\leq&\int_\Sigma e^{\f{bw^2}{\|w\|_{W^{1,2}(\Sigma,g)}^2}+\f{1}{b}{\|w\|_{W^{1,2}(\Sigma,g)}^2}}dv_g\\\nonumber
&\leq& C\int_\Sigma e^{\f{bw^2}{\|w\|_{W^{1,2}(\Sigma,g)}^2}}dv_g\\
&\leq& C.\label{T-C}
\eea
Thus (\ref{exp}) holds and the proof of Step 2 terminates.

 One can easily see that the conclusion of the lemma follows from (\ref{Harnack}). $\hfill\Box$\\

It is remarkable that only subcritical Trudinger-Moser inequality was employed in (\ref{T-C}).
Such inequalities are important tools in geometry and analysis. For
more details, we refer the reader to recent works \cite{Adi-Yang,Li-Yang,Nguyen,YangJDE,YangJGA,Yang-ZhuJFA,CR,Mancini} and the references therein.

\subsection{Palais-Smale condition}

For any $\lambda\in\mathbb{R}$, we define a functional $E_\lambda:W^{1,2}(\Sigma,g)\ra \mathbb{R}$ by
\be\label{functional}E_\lambda(u)=\int_\Sigma(|\nabla_gu|^2+2\kappa u-K_\lambda e^{2u})dv_g,\ee
where $\kappa:\Sigma\ra\mathbb{R}$ is the Gaussian curvature of $g$, $K_\lambda=K+\lambda$ is defined as in Theorem \ref{mainTheorem}.

\begin{lemma}[Palais-Smale condition]\label{PS}
Suppose that  $\Sigma_\lambda^-=\{x\in\Sigma: K_\lambda<0\}$ is nonempty for some $\lambda\in\mathbb{R}$. Then
$E_\lambda$ satisfies the $(PS)_c$ condition for all $c\in\mathbb{R}$, i.e., if $u_j$ is a sequence of functions in
$W^{1,2}(\Sigma,g)$ such that $E_\lambda(u_j)\ra c$ and $dE_\lambda(u_j)\ra 0$, then there exists some $u_0\in W^{1,2}(\Sigma,g)$
satisfying $u_j\ra u_0$ in $W^{1,2}(\Sigma,g)$.
\end{lemma}

\noindent{\it Proof.} Let $(u_j)$ be a function sequence such that $E_\lambda(u_j)\ra c$ and $dE_\lambda(u_j)\ra 0$, or equivalently
\bea\label{-c}&&\int_\Sigma(|\nabla_gu_j|^2+2\kappa u_j-K_\lambda e^{2u_j})dv_g=c+o_j(1),\\
\label{prim-0}&&\int_\Sigma(\nabla_gu_j\nabla_g\varphi+\kappa\varphi-K_\lambda e^{2u_j}\varphi)dv_g=o_j(1)\|\varphi\|_{W^{1,2}(\Sigma,g)},
\quad\forall \varphi\in W^{1,2}(\Sigma,g),\eea
where $o_j(1)\ra 0$ as $j\ra\infty$.

Note that ${\rm supp}\,\bm\beta=\{p_1,\cdots,p_\ell\}$ is a set of finite points. $\Sigma_\lambda^-\setminus{\rm supp}\,\bm\beta$
must contain a domain $\Omega$ such that the closure of $\Omega$ is also contained in $\Sigma_\lambda^-\setminus{\rm supp}\,\bm\beta$.
In view of (\ref{g0}), there would exist two positive constants $C_1$ and $C_2$ depending only on $\Omega$ such that
$$C_1g_0\leq g\leq C_2g_0\quad{\rm on}\quad \Omega.$$
Denote $u_j^+=\max\{u_j,0\}$. Based on an argument of Ding and Liu (\cite{Ding-Liu}, Lemma 2),
where a mistake was corrected by Borer, Galimberti and Struwe (\cite{B-G-Stru}, Appendix),
for another domain $\Omega^\prime\subset\subset\Omega$, there exists a positive constant $C$ depending only on $C_1$, $C_2$ and
${\rm dist}_g(\Omega^\prime,\p\Omega)$ such that
\be\label{W-bd}\int_{\Omega^\prime}(|\nabla_gu_j^+|^2+{u_j^+}^2)dv_g\leq C.\ee

Taking $\varphi\equiv 1$ in (\ref{prim-0}), one has
$$\int_\Sigma K_\lambda e^{2u_j}dv_g-\int_\Sigma \kappa dv_g=o_j(1).$$
This together with the Gauss-Bonnet formula (\ref{Gauss-Bonnet}) gives
\be\label{int-k}\int_\Sigma K_\lambda e^{2u_j}dv_g=2\pi\chi(\Sigma,\bm\beta)+o_j(1).\ee
Inserting (\ref{int-k}) into (\ref{-c}), we conclude
\be\label{W12}\int_\Sigma(|\nabla_gu_j|^2+2\kappa u_j)dv_g=c+2\pi\chi(\Sigma,\bm\beta)+o_j(1).\ee

We now {\it claim} that $u_j$ is bounded in $L^{2}(\Sigma,g)$. Suppose not, there holds $\|u_j\|_{L^{2}(\Sigma,g)}\ra\infty$.
We set $v_j=u_j/\|u_j\|_{L^{2}(\Sigma,g)}$.
Note that
$$\int_\Sigma \kappa\f{u_j}{\|u_j\|^2_{L^{2}(\Sigma,g)}}dv_g=o_j(1).$$
This together with (\ref{W12}) leads to
\be\label{t-0}\int_\Sigma |\nabla_gv_j|^2dv_g=o_j(1).\ee
Hence $v_j$ is bounded in $W^{1,2}(\Sigma,g)$ and (\ref{t-0}) leads to
$v_j\ra \gamma$ in $W^{1,2}(\Sigma,g)$ for some constant $\gamma$. Since
$\|v_j\|_{L^2(\Sigma,g)}=1$, we have $\gamma\not=0$. It follows from (\ref{W12}) that
\be\label{000}\int_\Sigma \kappa v_jdv_g\leq o_j(1).\ee
Letting $j\ra\infty$ in (\ref{000}),
we obtain $2\pi\chi(\Sigma,\bm\beta)\gamma\leq 0$ by using the Gauss-Bonnet formula (\ref{Gauss-Bonnet}).
Since $\chi(\Sigma,\bm\beta)<0$ and $\gamma\not=0$, we have $\gamma>0$.
  On the other hand, we conclude by (\ref{W-bd}) that
$$\int_{\Omega^\prime}(|\nabla_gv_j^+|^2+{v_j^+}^2)dv_g=o_j(1),$$
which leads to $\gamma\leq 0$. This contradicts $\gamma>0$ and confirms our claim.

Since $u_j$ is bounded in $L^{2}(\Sigma,g)$, we have by (\ref{W12}) that $u_j$ is bounded
in $W^{1,2}(\Sigma,g)$. Up to a subsequence, we can assume $u_j$ converges to $u_0$  weakly in $W^{1,2}(\Sigma,g)$, strongly in
$L^s(\Sigma,g)$ for any $s>1$. A Trudinger-Moser inequality for surfaces with conical singularities \cite{Troyanov} implies that
$e^{2u_j}$ is bounded in $L^s(\Sigma,g)$ for any $s>1$. Hence  $e^{u_j}$ converges to $e^{u_0}$ in $L^s(\Sigma,g)$ for any $s>1$. This together with
 (\ref{prim-0}) leads to
 \bna\int_\Sigma|\nabla_gu_j|^2dv_g&=&\int_\Sigma (-\kappa u_0+K_\lambda e^{2u_0}u_0)dv_g+o_j(1)\\
 &=&\int_\Sigma|\nabla_gu_0|^2dv_g+o_j(1).\ena
This implies that $u_j\ra u_0$ in $W^{1,2}(\Sigma,g)$. $\hfill\Box$

\subsection{Upper and lower solutions principle}
Let $f:\Sigma\times\mathbb{R}\ra \mathbb{R}$ be a smooth function.
$u\in W^{1,2}(\Sigma,g)\cap C^2(\Sigma\setminus{\rm supp}\,\bm\beta)\cap C^0(\Sigma)$ is defined to be an upper (lower) solution
to the elliptic equation
\be\label{ellip}\Delta_gu+f(x,u)= 0,\ee
if $u$ satisfies
 $\Delta_gu+f(x,u)\geq (\leq)\, 0$
in the distributional sense on $\Sigma$ and  point-wisely in $\Sigma\setminus{\rm supp}\,\bm\beta$.

\begin{lemma}[Upper and lower solutions principle]\label{upper-lower}
 Suppose that $\psi,\varphi\in W^{1,2}(\Sigma,g)\cap C^2(\Sigma\setminus{\rm supp}\,\bm\beta)\cap C^0(\Sigma)$ are upper  and
 lower solutions to (\ref{ellip}) respectively, and that $\varphi\leq\psi$ on $\Sigma$. Then (\ref{ellip}) has a solution $u\in W^{1,2}(\Sigma,g)\cap C^2(\Sigma\setminus{\rm supp}\,\bm\beta)\cap C^0(\Sigma)$ with $\varphi\leq u\leq \psi$ on $\Sigma$.
\end{lemma}
\noindent{\it Proof.} We follow the lines of Kazdan and Warner \cite{KW1}.
Let $A$ be a constant such that $-A\leq\varphi\leq\psi\leq A$. Since $\Sigma$ is closed, one finds a sufficiently large
constant $c$ such that $G(x,t)=ct+f(x,t)$ is increasing in $t\in[-A,A]$ for any fixed $x\in\Sigma$. Define an elliptic operator
$Lu=\Delta_gu+cu$ for
$u\in W^{1,2}(\Sigma,g)\cap C^2(\Sigma\setminus{\rm supp}\,{\bm\beta})\cap C^0(\Sigma)$. Now we define
\bna
&&\varphi_0=\varphi,\,\,\,\varphi_j=L^{-1}(G(x,\varphi_{j-1})),\,\,\,\forall j\geq 1\\
&&\psi_0=\psi,\,\,\,\psi_j=L^{-1}(G(x,\psi_{j-1})),\,\,\,\forall j\geq 1.
\ena
Here $L^{-1}: L^2(\Sigma,g)\ra W^{1,2}(\Sigma,g)$ is well defined due to the Lax-Milgram theorem.
This together with the definition of upper and lower solutions and the monotonicity of $G(x,t)$ with respect to $t$
leads to
$$L\varphi\leq L\varphi_1=G(x,\varphi)\leq G(x,\psi)=L\psi_1\leq L\psi.$$
Then the weak maximum principle (Lemma \ref{maximum principle}) implies that
$$\varphi\leq\varphi_1\leq\psi_1\leq\psi.$$
By induction, we have
$$\varphi\leq\varphi_{j-1}\leq\varphi_j\leq\psi_j\leq\psi_{j-1}\leq\psi,\,\,\,j=1,2,\cdots.$$
Clearly we can assume that $\varphi_j$ converges to $u_1$ and $\psi_j$ converges to $u_2$ point-wisely.
By elliptic estimates, one concludes that the above convergence is in $C_{\rm loc}^2(\Sigma\setminus{\rm supp}\,\bm\beta)\cap
C^0(\Sigma)$. Moreover, $v=u_1$ or $u_2$ is a distributional solution to $Lv=G(x,v)$. $\hfill\Box$

\section{Proof of Theorem \ref{thm2}}

In this section, we prove Theorem \ref{thm2} by using variational method.

\subsection{Unique solution in the case $\lambda\leq 0$}

{\it Proof of $(i)$ of Theorem \ref{thm2}}. Assume $\max_\Sigma K=0$ and $K\not\equiv 0$.
If $\lambda<0$, this has been proved by Troyanov (\cite{Troyanov}, Theorem 1). We now consider the general case
$\lambda\leq 0$. Let $E_\lambda$ be the functional defined as in (\ref{functional}), where
$K_\lambda=K+\lambda$.\\

{\bf Claim 1}. {\it $E_\lambda$ is strict convex on $W^{1,2}(\Sigma,g)$}.\\

It suffices to prove that for any $u\in W^{1,2}(\Sigma)$, there exists some constant $C>0$  such that
\be\label{convex}d^2E_\lambda(u)(h,h)\geq C\|h\|_{W^{1,2}(\Sigma,g)}^2\quad\forall h\in W^{1,2}(\Sigma,g).\ee
Suppose not. There would be a function $u\in W^{1,2}(\Sigma,g)$ and
a function sequence $(h_j)\subset W^{1,2}(\Sigma,g)$ such that $\|h_j\|_{W^{1,2}(\Sigma,g)}=1$ for all $j$ and
$d^2E_\lambda(u)(h_j,h_j)\ra 0$ as $j\ra\infty$.
One may assume up to a subsequence, $h_j$ converges to $h_\infty$ weakly in $W^{1,2}(\Sigma,g)$, strongly in $L^p(\Sigma,g)$ for any
$p>1$, and almost everywhere in $\Sigma$. Since
$$d^2E_\lambda(u)(h_j,h_j)=2\int_\Sigma(|\nabla_gh_j|^2-2K_\lambda e^{2u}h_j^2)dv_g$$ and $K_\lambda\leq 0$, we conclude
$\int_\Sigma|\nabla_gh_j|^2dv_g\ra 0$ and $\int_\Sigma K_\lambda e^{2u}h_j^2dv_g\ra 0$, which leads to
$h_\infty\equiv C_0$ for some constant $C_0$, and further
$$\label{c0}C_0^2\int_\Sigma K_\lambda e^{2u}dv_g=\int_\Sigma K_\lambda e^{2u}h_\infty^2dv_g=\lim_{j\ra\infty}\int_\Sigma
K_\lambda e^{2u}h_j^2dv_g=0.$$
Clearly $\int_\Sigma K_\lambda e^{2u}dv_g<0$, and thus $C_0=0$. This contradicts
$\|h_\infty\|_{L^2(\Sigma,g)}=\lim_{j\ra\infty}\|h_j\|_{L^{2}(\Sigma,g)}= 1$. Hence (\ref{convex}) holds.\\

{\bf Claim 2}. {\it $E_\lambda$ is coercive}.\\

Since for any $\epsilon>0$, there exists a constant $C(\epsilon)$ such that
$\int_\Sigma\kappa udv_g\leq \epsilon\|u\|_{W^{1,2}(\Sigma,g)}^2+C(\epsilon)$, it suffices to find some constant $C>0$ such that for all $u\in W^{1,2}(\Sigma,g)$, there holds
\be\label{geq}\int_\Sigma(|\nabla_gu|^2-K_\lambda e^{2u})dv_g\geq C\|u\|_{W^{1,2}(\Sigma,g)}^2.\ee
Suppose not. There would exist a sequence of functions $(u_j)$ satisfying
$$\int_\Sigma(|\nabla_g u_j|^2+u_j^2)dv_g=1,\quad\int_\Sigma(|\nabla_gu_j|^2-K_\lambda e^{2u_j})dv_g=o_j(1).$$ It follows that
up to a subsequence,
$u_j$ converges to $u^\ast$ weakly in $W^{1,2}(\Sigma,g)$ and strongly in $L^p(\Sigma,g)$ for any $p>1$. One easily see that
 $$0< \int_\Sigma(|\nabla_gu^\ast|^2-K_\lambda e^{2u^\ast})dv_g\leq \lim_{j\ra\infty}
 \int_\Sigma(|\nabla_gu_j|^2-K_\lambda e^{2u_j})dv_g=0,$$ which is impossible. Hence (\ref{geq}) holds.

 In view of {\it Claims} 1 and 2, a direct method of variation shows $\inf_{u\in W^{1,2}(\Sigma,g)}E_\lambda(u)$ can be attained
 by some $u_0\in W^{1,2}(\Sigma,g)$ and $u_0$ is
 the unique critical point of $E_\lambda$. $\hfill\Box$

 \subsection{Existence of $\lambda^\ast$}
 When $\lambda=0$, the equation (\ref{eqn}) becomes
 \be\label{lm-0}\Delta_gu+\kappa-Ke^{2u}=0\quad {\rm on}\quad \Sigma.\ee
 Let $u$ be a solution of (\ref{lm-0}). The linearized equation of (\ref{lm-0}) at $u$ reads
 $\Delta_gv-2K e^{2u}v=0$, which has a unique solution $v\equiv 0$. By the implicit theorem, there is a sufficiently
 small $s>0$ such that for any $\lambda\in(0,s)$, the equation (\ref{eqn}) has a solution. Define
 \be\label{lm-ast}\lambda^\ast=\sup\le\{s:\, {\rm the\, equation}\, (\ref{eqn})
 \,{\rm has\,a\,solution\,for\,any}\, \lambda\in(0,s)\ri\}.\ee
 One can see that $\lambda^\ast\leq -\min_\Sigma K$. For otherwise $K_\lambda> 0$ for some $\lambda<\lambda^\ast$.
 Integrating (\ref{eqn}), we obtain
 $$0>2\pi\chi(\Sigma,\bm{\beta})=\int_\Sigma \kappa dv_g=\int_\Sigma K_\lambda e^{2u}dv_g\geq 0,$$
 which is impossible. In conclusion, we have
 $0<\lambda^\ast\leq-\min_\Sigma K$. Further analysis (Subsection 4.4, Claim 2) implies that
 $\lambda^\ast<-\min_\Sigma K$.

 \subsection{Multiplicity of solutions for $0<\lambda<\lambda^\ast$\\}

  {\it Proof of $(ii)$ of Theorem \ref{thm2}}. Fix $\lambda$, $0<\lambda<\lambda^\ast$. We shall seek two
  different solutions of
  (\ref{eqn}), one is a strict local minimum of the functional $E_\lambda$, the other is of the mountain-pass type.
  The proof will be divided into several steps below.\\

  {\bf Step 1}. {\it Existence of upper and lower solutions}.\\

  Take $\lambda_1$ with $\lambda<\lambda_1<\lambda^\ast$. Let $u_{\lambda_1}\in
  W^{1,2}(\Sigma,g)\cap C^2(\Sigma\setminus{\rm supp}\,\bm\beta)\cap C^0(\Sigma)$ be a solution of (\ref{eqn}) at $\lambda_1$.
  Set $\psi=u_{\lambda_1}$. One can see that $\psi$ is a strict upper solution of (\ref{eqn}), namely
  \be\label{upp}\Delta_g\psi+\kappa-K_{\lambda}e^{2\psi}>0.\ee
  Clearly the equation
  \be\label{eta}\Delta_g\eta=-\kappa+\f{1}{{\rm Vol}_g(\Sigma)}\int_\Sigma \kappa dv_g\ee
  has a distributional solution $\eta\in W^{1,2}(\Sigma,g)\cap C^2(\Sigma\setminus{\rm supp}\,\bm\beta)\cap C^0(\Sigma)$.
  Let $\varphi=\eta-s$, where $s$ is a positive constant.
 Obviously $\varphi< \psi$ on $\Sigma$ for sufficiently large $s$. Since $\int_\Sigma \kappa dv_g=2\pi\chi(\Sigma,\bm{\beta})<0$, we have
 \be\label{lower}\Delta_g\varphi+\kappa-K_\lambda e^{2\varphi}=\f{1}{{\rm Vol}_g(\Sigma)}\int_\Sigma \kappa dv_g-K_\lambda e^{2\eta-2s}<0,\ee
 provided that $s$ is chosen sufficiently large. Thus $\varphi$ is a strict lower solution of (\ref{eqn}).\\

 {\bf Step 2}. {\it The first solution of (\ref{eqn}) can be chosen as a strict local minimum of $E_\lambda$.}\\

 Let $f_\lambda(x,t)=ct-\kappa+K_\lambda e^{2t}$. Fix  a sufficiently large positive constant $c$
 such that $f_\lambda(x,t)$ is increasing in $t\in[-A,A]$, where $A$ is a constant such that $-A\leq\varphi<\psi\leq A$. Let
 $F_\lambda(x,u)=\int_0^uf_\lambda(x,t)dt$.
 It is easy to see that
 $$E_\lambda(u)=\int_\Sigma|\nabla_gu|^2dv_g+c\int_\Sigma u^2dv_g-2\int_\Sigma F_\lambda(x,u)dv_g-\int_\Sigma K_\lambda dv_g.$$
  %{\it Proof.} Let $$\mathscr{S}=\le\{u\in W^{1,2}(\Sigma)\cap C^2(\Sigma\setminus
 %\{p_1,\cdots,p_\ell\})\cap C^0(\Sigma): \varphi(x)\leq \psi(x),\forall x\in\Sigma\ri\}.$$
 %{\bf (Maybe we should write $\bm{W^{2,p}(\Sigma)}$ for some $\bm {p>1}$ instead of
 %$\bm{C^2(\Sigma\setminus
 %\{p_1,\cdots,p_\ell\})\cap C^0(\Sigma)}$!)}
 Define a function
 $$\hat{f}_\lambda(x,t)=\le\{\begin{array}{lll}
 f_\lambda(x,\psi(x))&{\rm when}&t>\psi(x)\\[1.5ex]
 f_\lambda(x,t)&{\rm when}&\varphi(x)\leq t\leq\psi(x)\\[1.5ex]
 f_\lambda(x,\varphi(x))&{\rm when}&t<\varphi(x)
 \end{array}\ri.$$
 and a functional
 $$\hat{E}_\lambda(u)=\int_\Sigma(|\nabla_gu|^2+cu^2)dv_g-2\int_\Sigma \hat{F}_\lambda(x,u)dv_g-\int_\Sigma K_\lambda dv_g,$$
 where $\hat{F}_\lambda(x,t)=\int_0^t\hat{f}_\lambda(x,s)ds$.
 Obviously $\hat{E}_\lambda$ is bounded from below on $W^{1,2}(\Sigma,g)$. Denote
 $$a=\inf_{u\in W^{1,2}(\Sigma,g)}\hat{E}_\lambda(u).$$
 Taking a function sequence $(u_j)\subset W^{1,2}(\Sigma,g)$ such that $\hat{E}_\lambda(u_j)\ra a$ as $j\ra\infty$.
 It follows that $u_j$ is bounded in $W^{1,2}(\Sigma,g)$, and thus up to a subesequence the Sobolev embedding and the Trudinger-Moser
 inequality
 lead to $u_j$ converges to some $u_\lambda$ weakly in $W^{1,2}(\Sigma,g)$, strongly in $L^q(\Sigma,g)$ for any $q>1$, almost
 everywhere in $\Sigma$, and
 $e^{2u_j}$ converges to $e^{2u_\lambda}$ in $L^1(\Sigma,g)$. Hence $\hat{E}_\lambda(u_\lambda)\leq a$. Then by the definition of $a$, we
 conclude
 $$\hat{E}_\lambda(u_\lambda)=\inf_{u\in W^{1,2}(\Sigma,g)}\hat{E}_\lambda(u).$$
 As a consequence $u_\lambda$ satisfies the Euler-Lagrange equation
 \be\label{u-l}
 \Delta_g u_\lambda+cu_\lambda=\hat{f}_\lambda(x,u_\lambda)\ee
 in the distributional sense. By elliptic estimates, one has $u_\lambda\in C^2(\Sigma\setminus
 {\rm supp}\,\bm\beta)\cap C^0(\Sigma)$.

 Noting that $f(x,t)$ is increasing with respect to $t\in [-A,A]$, we have
 $$\Delta_g \varphi(x)+c\varphi(x) \leq f_\lambda(x,\varphi(x))\leq \hat{f}_\lambda(x,u_\lambda(x))\leq f_\lambda(x,\psi(x))\leq
 \Delta_g \psi(x)+c\psi(x)$$
 in the distributional sense. In view of (\ref{upp}), (\ref{lower}) and (\ref{u-l}), one concludes by the strong maximum principle
 (Lemma \ref{strong-maximum}) that
 \be\label{inner-point}\varphi(x)<u_\lambda(x)<\psi(x),\quad\forall x\in\Sigma.\ee
 Obviously $\hat{E}_\lambda(u)=E_\lambda(u)$ for all $u\in W^{1,2}(\Sigma)$ with $\varphi\leq u\leq \psi$.
 For any $h\in C^1(\Sigma)$, we define a function
 $\zeta(t)=E(u_\lambda+th)$, $t\in\mathbb{R}$. In view of (\ref{inner-point}), there holds $\varphi\leq u_\lambda+th\leq\psi$
 and thus $\hat{E}_\lambda(u_\lambda+th)=E_\lambda(u_\lambda+th)$, provided that $|t|$ is sufficiently small.
 Since $u_\lambda$ is a minimum of $\hat{E}_\lambda$ on $W^{1,2}(\Sigma,g)$, we have $\zeta^\prime(0)=dE_\lambda(u_\lambda)(h)=0$ and $\zeta^{\prime\prime}(0)=d^2E_\lambda(u_\lambda)(h,h)\geq 0$.
 Therefore we have
 \bea\label{crit}
 \int_\Sigma(\nabla_gu_\lambda\nabla_gh+\kappa h-K_\lambda e^{2u_\lambda}h)dv_g=0,&\forall h\in C^1(\Sigma),&\\[1.2ex]
 \label{2-order}\int_\Sigma(|\nabla_gh|^2-2K_\lambda e^{2u_\lambda}h^2)dv_g\geq 0,&\forall h\in C^1(\Sigma).&
 \eea
 Since $C^1(\Sigma)$ is dense in $W^{1,2}(\Sigma,g)$, (\ref{crit}) and (\ref{2-order}) still hold for all $h\in W^{1,2}(\Sigma,g)$.
 We further prove that
 there exists a positive constant $C$ such that
 \be\label{strict}d^2{E}_\lambda(u_\lambda)(h,h)\geq C\|h\|^2_{W^{1,2}(\Sigma,g)},\quad\forall h\in W^{1,2}(\Sigma,g).\ee
  For the proof of (\ref{strict}), we adapt an argument of Borer, Galimberti, and Struwe (\cite{B-G-Stru}, Section 2).
  Since
 $d^2{E}_\lambda(u_\lambda)(h,h)\geq 0$
 for all $h\in W^{1,2}(\Sigma,g)$, we have
 $$\Lambda:=\inf_{\|h\|_{W^{1,2}(\Sigma,g)}=1}d^2{E}_\lambda(u_\lambda)(h,h)\geq 0.$$
 Suppose $\Lambda=0$. We {\it claim} that  there exists some $h$ with $\|h\|_{W^{1,2}(\Sigma,g)}=1$ such that
 $d^2{E}_\lambda(u_\lambda)(h,h)=0$. To see this, we let $h_j$ satisfy $\|h_j\|_{W^{1,2}(\Sigma,g)}=1$ and
 $d^2{E}_\lambda(u_\lambda)(h_j,h_j)\ra 0$ as $j\ra\infty$. Up to a subsequence, we can assume $h_j$ converges to some $h$
 weakly in $W^{1,2}(\Sigma,g)$, strongly in $L^q(\Sigma,g)$ for all $q>1$, and almost everywhere in $\Sigma$.
 It follows that
 $$\lim_{j\ra\infty}\int_\Sigma|\nabla_gh_j|^2dv_g=\int_\Sigma 2K_\lambda e^{2u_\lambda}h^2dv_g\leq \int_\Sigma|\nabla_gh|^2dv_g.$$
 This leads to $h_j\ra h$ in $W^{1,2}(\Sigma,g)$ as $j\ra\infty$, and confirms our claim.
   Moreover, since the functional $v\mapsto d^2{E}_\lambda(u_\lambda)(v,v)$ attains its minimum at $v=h$,
 it follows that
 $d^2E_\lambda(u_\lambda)(h,w)=0$ for all $w\in W^{1,2}(\Sigma,g)$;
 that is, $h$ is a weak solution of the equation
 \be\label{22}\Delta_gh=2K_\lambda e^{2u_\lambda}h.\ee
 Note that $h$ is not a constant. For otherwise (\ref{22}) yields
 $$0>2\pi\chi(\Sigma,\bm\beta)=\int_\Sigma K_\lambda e^{2u_\lambda}dv_g=0,$$
 which is impossible. Multiplying (\ref{22}) by $h^3$, we get
 $$d^4E_\lambda(u_\lambda)(h,h,h,h)=-16\int_\Sigma K_\lambda e^{2u_\lambda}h^4dv_g=-24
 \int_\Sigma h^2|\nabla_gh|^2dv_g<0.$$
 Since $d^2E_\lambda(u_\lambda+th)(h,h)$ attains its minimum at $t=0$, we have
 $d^3E_\lambda(u_\lambda)(h,h,h)=0$, which together with the facts $dE_\lambda(u_\lambda)=0$ and $d^2E_\lambda(u_\lambda)(h,h)=0$
 leads to
 \be\label{O-e}E_\lambda(u_\lambda+\epsilon h)=E_\lambda(u_\lambda)+\f{\epsilon^4}{24}d^4E_\lambda(u_\lambda)(h,h,h,h)
 +O(\epsilon^5)<E_\lambda(u_\lambda)\ee
 for small $\epsilon>0$. Applying elliptic estimates to (\ref{22}), we have $h\in C^0(\Sigma)$. Then there exists $\epsilon_0>0$
  such that if $0<\epsilon<\epsilon_0$, then $\varphi\leq u_\lambda+\epsilon h\leq \psi$ on $\Sigma$, and thus by (\ref{O-e}),
  $$\hat{E}_\lambda(u_\lambda+\epsilon h)=E(u_\lambda+\epsilon h)<E_\lambda(u_\lambda)=\hat{E}_\lambda(u_\lambda),$$
  contradicting the fact that $u_\lambda$ is the minimum of $\hat{E}_\lambda$.
 Therefore $\Lambda>0$ and (\ref{strict}) follows immediately. As a consequence, $u_\lambda$ is a strict local minimum of $E_\lambda$ on $W^{1,2}(\Sigma,g)$.\\

 {\bf Step 3}. {\it The second solution of (\ref{eqn}) can be achieved by a mountain pass theorem.}\\

 Let $u_\lambda$ be as in Step 2. Since $u_\lambda$ is a strict local minimum of $E_\lambda$ on $W^{1,2}(\Sigma,g)$,
  there would exist a sufficiently small $r>0$ such that
    \be\label{g-e}\inf_{\|u-u_\lambda\|_{W^{1,2}(\Sigma,g)}=r}E_\lambda(u)>E_\lambda(u_\lambda).\ee
 Moreover, a calculation of Ding and Liu (\cite{Ding-Liu}, Page 1061) shows
 for any $\lambda>0$, $E_\lambda$ has no lower bound on $W^{1,2}(\Sigma,g)$.
 In particular, there exists some $v\in W^{1,2}(\Sigma,g)$ verifying that
 \be\label{e}E_\lambda(v)<E_\lambda(u_\lambda),\quad \|v-u_\lambda\|_{W^{1,2}(\Sigma,g)}>r.\ee
 Combining (\ref{g-e}), (\ref{e}) and Lemma \ref{PS}, we obtain by using the mountain-pass theorem due to
 Ambrosetti and Rabinowitz \cite{Ambrosetti-Rabinowitz} that
 the mini-max value $$c=\min_{\gamma\in\Gamma}\max_{u\in\gamma}E_\lambda(u)$$
 is a critical value of $E_\lambda$, where
 $\Gamma=\{\gamma\in \mathcal{C}([0,1],W^{1,2}(\Sigma,g)):\gamma(0)=u_\lambda,\gamma(1)=v\}$.
 Equivalently there exists some $u^\lambda\in W^{1,2}(\Sigma,g)$ satisfying
 $E_\lambda(u^\lambda)=c$ and $dE_\lambda(u^\lambda)=0$. Thus $u^\lambda$ is a solution of the equation (\ref{eqn})
 and $u^\lambda\not=u_\lambda$.
  Finally elliptic estimates imply that $u^\lambda\in C^2(\Sigma\setminus{\rm supp}\,\bm\beta)\cap C^0(\Sigma)$.
  $\hfill\Box$

 \subsection{Solvability of (\ref{eqn}) at $\lambda^\ast$}

 {\it Proof of $(iii)$ of Theorem \ref{thm2}.}
For any $\lambda$, $0<\lambda<\lambda^\ast$, we let $u_\lambda$ be the local minimum of $E_\lambda$ obtained in the previous
subsection. In particular, $u_\lambda$ is a solution of
(\ref{eqn}) and
$$\int_\Sigma(|\nabla_g\phi|^2-2K_\lambda e^{2u_\lambda}\phi^2)dv_g\geq 0,
\quad \forall \phi\in W^{1,2}(\Sigma,g).$$
The remaining part of the proof will be divided into several claims as below.\\

{\bf Claim 1}. {\it There exists some constant $C$ such that $u_\lambda\geq -C$ on $\Sigma$ uniformly in $\lambda\in(0,\lambda^\ast)$.}\\

To see this, we let $\eta$ satisfy (\ref{eta}) and $\varphi_{s}=\eta-s$ for $s>0$. The analog of (\ref{lower}) reads
$$\Delta_g\varphi_{s}+\kappa-Ke^{2\varphi_s}<0$$
with $s$ chosen sufficiently large, say $s\geq s_0$. Equivalently $\varphi_s$ is a lower solution of (\ref{eqn}) at $\lambda=0$,
provided that $s\geq s_0$. Clearly $\varphi_{s}$ is also a strict lower solution of (\ref{eqn}) at $\lambda\in(0,\lambda^\ast)$ for
any $s\geq s_0$.
 We now prove that $u_\lambda\geq\varphi_{s_0}$, and consequently claim 1 holds. For otherwise, by varying $s\in[s_0,\infty)$, we
 find that for some $s$ there holds $u_\lambda\geq \varphi_s$ on $\Sigma$, and $u_\lambda(x_0)=\varphi_s(x_0)$ for some $x_0\in \Sigma$.
 Then the strong maximum principle (Lemma \ref{strong-maximum}) implies that $u_\lambda\equiv \varphi_s$ on $\Sigma$, which is impossible.\\

 {\bf Claim 2}. {\it Let $\Sigma_{\lambda^\ast}^-=\{x\in\Sigma: K_{\lambda^\ast}(x)<0\}$.
 Then $\Sigma_{\lambda^\ast}^-\not=\varnothing$.}\\

 Suppose $K_{\lambda^\ast}\geq 0$.  Let $\widetilde{g}=e^{2v}g$ be a metric with constant Gaussian curvature $-1$,
 where ${v}$ is a solution of $\Delta_g{v}+\kappa+e^{2{v}}=0$.
  In view of $(i)$ of Theorem \ref{thm2}, such a function $v$ uniquely exists.
 Let $w_\lambda=u_\lambda-v$. Noting that $\Delta_g=e^{2v}\Delta_{\widetilde{g}}$, we have
  $$\Delta_{\widetilde{g}}w_\lambda-1-K_\lambda e^{2w_\lambda}=0.$$
  Multiplying the above equation by $e^{-2w_\lambda}$ and integrating by parts, one has
  $$\int_\Sigma K_\lambda dv_{\widetilde{g}}=-\int_\Sigma e^{-2w_\lambda}dv_{\widetilde{g}}-2\int_\Sigma
  |\nabla_{\widetilde{g}}w_\lambda|^2e^{-2w_\lambda}dv_{\widetilde{g}}\leq 0.$$
  Hence
  $$\int_\Sigma K_{\lambda^\ast}dv_{\widetilde{g}}=\lim_{\lambda\ra\lambda^\ast}\int_\Sigma K_{\lambda}dv_{\widetilde{g}}\leq 0.$$
  This together with $K_{\lambda^\ast}\geq 0$ leads to $K_{\lambda^\ast}\equiv 0$, which contradicts the assumption that
  $K_{\lambda^\ast}$ is not a constant.\\

  {\bf Claim 3}. {\it Let $\Omega$ and $\Omega^\prime$ are two domains in $\Sigma$ such that
  $\Omega^\prime\subset\subset\Omega\subset\subset \Sigma_{\lambda^\ast}^-\setminus{\rm supp}\,\bm\beta$. Then
  $u_\lambda^+$ is bounded in $W^{1,2}(\Omega^\prime,g)$ with respect to $\lambda\in(0,\lambda^\ast)$.}\\

  Note that $K: \Sigma\ra\mathbb{R}$ is H\"older continuous. If $\lambda\in(0,\lambda^\ast)$, then
  $$\sup_{\Omega}K_\lambda\leq\sup_{\Omega}K_{\lambda^\ast}\leq -\epsilon$$
  for some $\epsilon>0$ depending only on $K$, $\lambda^\ast$ and $\Omega$.
  Similar to the proof of (\ref{W-bd}), we conclude Claim 3.\\

  {\bf Claim 4}. {\it The equation (\ref{eqn}) is solvable at $\lambda^\ast$.} \\

  Having Claims 1-3 in hand and arguing as Ding and Liu did in the proof of (\cite{Ding-Liu}, $(c)$ of the main theorem),
  we conclude that both $e^{2u_\lambda}$ and $u_\lambda$ are
   bounded in $L^q(\Sigma,g)$ for all $q>1$. By elliptic estimates, we have up to a subsequence, $u_\lambda$ converges to
   some $u$ in $W^{1,2}(\Sigma,g)$, where $u$ is a solution of
   $$\Delta_gu+\kappa-K_{\lambda^\ast}e^{2u}=0.$$
   By elliptic estimates, $u\in C^2(\Sigma\setminus{\rm supp}\,\bm\beta)\cap C^0(\Sigma)$. This gives the desired result. $\hfill\Box$

 \subsection{The equation (\ref{eqn}) has no distributional solution when $\lambda>\lambda^\ast$}

 {\it Proof of $(iv)$ of Theorem \ref{thm2}.} Suppose (\ref{eqn}) has a solution $u_{\lambda_1}$ at some $\lambda_1>\lambda^\ast$.
   Then for any ${\lambda}, 0<\lambda<\lambda_1$, $u_{\lambda_1}$
   is an upper solution of (\ref{eqn}). Similar to (\ref{lower}), we can easily construct a lower solution $\varphi$ of (\ref{eqn})
   such that $\varphi\leq u_{\lambda_1}$. In view of the upper and lower solutions principle (Lemma \ref{upper-lower}),
   there would exist a solution of (\ref{eqn}), which contradicts
    the definition of $\lambda^\ast$ (see (\ref{lm-ast}) above). $\hfill\Box$\\

{\bf Acknowledgements}. This work is supported by National Science Foundation of China (Grant Nos. 11171347,
 11471014, 41275063 and 11401575).

\end{document}